\documentclass{amsart}
\usepackage{amsfonts}
\usepackage{cases}

\RequirePackage{amsmath} \RequirePackage{amssymb}
\usepackage{amscd,latexsym,amsthm,amsfonts,amssymb,amsmath,amsxtra}
\usepackage
{hyperref}

    \theoremstyle{plain}
    \newtheorem{thm}{Theorem}[section] \newtheorem{cor}[thm]{Corollary}
    \newtheorem{lem}[thm]{Lemma}

    \numberwithin{equation}{section}
\title{On the Asymptotic Formula of $L'(1,\chi)$}
\author{Luhao Yan}
\begin{document}
\maketitle
\begin{abstract}
Let $\chi$ be a quadratic Dirichlet character. In some literatures, various asymptotic formulae of $L'(1,\chi)$,
under the assumption that $L(1,\chi)$ takes a small value, were derived. In this paper, we will give a new treatment unified for the odd and even cases, not depending on Kronecker limit formula. For imaginary quadratic fields, our result coincides with Proposition 22.10 in \cite{IK}.
\end{abstract}

\section{Introduction}
Siegel zeros (if possible) cause many strange phenomena. For a good survey, we refer to \cite{I}.
In \cite{G},\cite{GS} and Chapter 22 of \cite{IK}, some asymptotic formulae were derived. And for imaginary quadratic fields, our following main result is the same as Proposition 22.10 in \cite{IK} essentially.
\begin{thm}\label{main4}
Let $\chi$ be a primitive quadratic Dirichlet character modulo $q$.
If $L(1,\chi)\ll (\log q)^{-26}$, then we have
$$L'(1,\chi)=\frac{\pi^2}{6}\prod_{p|q}\left(1+\frac{1}{p}\right)\prod_{p\leqslant q\atop \chi(p)=1}\left(1+\frac{1}{p}\right)\left(1-\frac{1}{p}\right)^{-1}\cdot\Big(1+\mathcal{O}\left((\log q)^{-1/10}\right)\Big).$$
\end{thm}

\section{The Proof of the Main Theorem}
We need some preliminary results before proving it.
\begin{lem}\label{lem:ide}
Suppose $f(m,n)$ is an arithmetic function in two variables, and $x>u\geqslant 1$. Then
\begin{eqnarray*}
  \sum_{k<\sqrt{\frac{x}{u}}}\sum_{u<n\leqslant \frac{x}{k^2}}f(k^2,n) &=& \sum_{d\leqslant u}\sum_{u<n\leqslant \frac{x}{d}}\sum_{r\leqslant \frac{x}{dn}}\lambda(d)f(dr,n) \\
    & & {}+\sum_{u<n\leqslant x}\sum_{u<m\leqslant \frac{x}{n}}\sum_{d\mid m\atop d>u}\lambda(d)f(m,n),
\end{eqnarray*}
where $\lambda(n)=(-1)^{\Omega(n)}$~.
\end{lem}
Proof: Using the following property of the function $\lambda(n)$:
\begin{equation*}
  \sum_{d\mid m}\lambda(d)=
  \begin{cases}
     1, & \text{$m$ is a square,}\\
     0, & \text{otherwise.}
   \end{cases}
\end{equation*}
we get that
\begin{eqnarray*}
  \sum_{k<\sqrt{\frac{x}{u}}}\sum_{u<n\leqslant \frac{x}{k^2}}f(k^2,n) &=& \sum_{u<mn\leqslant x\atop u<n}f(m,n)\sum_{d\mid m}\lambda(d) \\
    &=& \sum_{u<mn\leqslant x\atop u<n}f(m,n)\sum_{d\mid m\atop d\leqslant u}\lambda(d)+\sum_{u<mn\leqslant x\atop u<n}f(m,n)\sum_{d\mid m\atop d>u}\lambda(d) \\
    &=& \sum_{d\leqslant u}\sum_{u<n\leqslant \frac{x}{d}}\sum_{r\leqslant \frac{x}{dn}}\lambda(d)f(dr,n)\\
    & &{}+\sum_{u<n\leqslant x}\sum_{u<m\leqslant \frac{x}{n}}\sum_{d\mid m\atop d>u}\lambda(d)f(m,n).
\end{eqnarray*}

\begin{lem}\label{lem:psi}
Suppose $\chi$ is a primitive real Dirichlet character modulo $q$ and $q<x$.
Let $\psi_u(z,\chi)=\sum\limits_{u<n\leqslant z}\Lambda(n)\chi(n)$, for $z>u$.
we have
\begin{eqnarray*}
  \sum_{k<\sqrt{\frac{x}{u}}\atop (k,q)=1}\psi_u\left(\frac{x}{k^2},\chi\right) &=& \sum_{u<m\leqslant x/u}\rho_u(m)\chi(m)\psi_u\left(\frac{x}{m},\chi\right) \\
    & & {}+\mathcal{O}\left(\left(q\sqrt q+\frac{x}{\sqrt q}+u^2\sqrt q\right)\log^2 x\right),
\end{eqnarray*}
where $$\rho_u(m)=\sum_{d\mid m\atop d>u}\lambda(d).$$
\end{lem}

Proof:
Let $\alpha=a/q$ with $(a,q)=1$, and $e(t)=\exp(2\pi it)$.\\
Set $f(m,n)=\Lambda(n)e(\alpha mn).$ By Lemma \ref{lem:ide}, we have

\begin{eqnarray}\label{rep}
\sum_{k<\sqrt{\frac{x}{u}}}\sum_{u<n\leqslant \frac{x}{k^2}}\Lambda(n)e(\alpha k^2 n) &=& \sum_{d\leqslant u}\sum_{u<n\leqslant \frac{x}{d}}\sum_{r\leqslant \frac{x}{dn}}\lambda(d)\Lambda(n)e(\alpha drn)\\
    & & {}+\sum_{u<n\leqslant x}\sum_{u<m\leqslant \frac{x}{n}}\rho_u(m)\Lambda(n)e(\alpha mn)\nonumber\\
    &:=& T^\sharp(\alpha)+T^\flat(\alpha).\nonumber
\end{eqnarray}

Let $rn=l$,
\begin{eqnarray*}
T^\sharp(\alpha)&=& \sum_{d\leqslant u}\lambda(d)\sum_{l\leqslant x/d}e(\alpha
dl)\sum_{n>u,\;n|l}\Lambda(n)\\
&=& \sum_{d\leqslant u}\lambda(d)\sum_{l\leqslant x/d}e(\alpha dl)\log l-\sum_{d\leqslant u}\lambda(d)\sum_{l\leqslant x/d}e(\alpha dl)\sum_{n\leqslant u,\;n|l}\Lambda(n)\\
&:=& T^\sharp_1(\alpha)-T^\sharp_2(\alpha).
\end{eqnarray*}

Applying the basic estimations $$\left|\sum\limits_{1\leq n\leq N}e(\alpha
n)\right|\leqslant \min \left(N,\frac {1}{2\|\alpha\|}\right)$$
and
$$\sum\limits_{1\leq n\leq N}\min \left(\frac xn,\frac {1}{2\|\alpha
n\|}\right)\ll N\log q+\frac xq\log N+q\log q,$$
we derive by partial summation that
\begin{eqnarray*}
T^\sharp_1(\alpha)&\ll&\sum_{d\leqslant u}\left|\sum_{l\leqslant x/d}e(\alpha
dl)\log l\right|\\
 &\ll&\log x\sum_{d\leqslant u}\min \left(\frac xd,\;\frac
{1}{2\|\alpha d\|}\right)\\
&\ll&(xq^{-1}+u+q)\log^2 x.
\end{eqnarray*}

For $T^\sharp_2(\alpha)$, let $l=rn$, we get
\begin{eqnarray*}
T^\sharp_2(\alpha)&\ll&\sum_{d\leqslant u}\sum_{n\leqslant u}\Lambda(n)
\left|\sum_{r\leqslant \frac{x}{dn}}e(\alpha drn)\right|\\
&\ll&\sum_{d\leqslant u}\sum_{n\leqslant u}\Lambda(n)\min \left(\frac {x}{dn},\;\frac
{1}{2\|\alpha dn\|}\right)\\
&\ll&\sum_{h\leqslant u^2}\min \left(\frac {x}{h},\frac {1}{2\|\alpha h\|}\right)\sum_{n|h}\Lambda(n)\\
&\ll&(xq^{-1}+u^2+q)\log^2 x.
\end{eqnarray*}

For a real primitive character $\chi$ modulo $q$,
\begin{eqnarray*}
\sum_{k<\sqrt{\frac{x}{u}}\atop (k,q)=1}\psi_u\left(\frac{x}{k^2},\chi\right)&=& \sum_{k<\sqrt{\frac{x}{u}}}\sum_{u<n\leqslant \frac{x}{k^2}}
\Lambda(n)\chi(k^2n)\\
&=& \frac {1}{\tau(\chi)}\sum_{a=1}^{q}\chi(a)\sum_{k<\sqrt{\frac{x}{u}}}\sum_{u<n\leqslant \frac{x}{k^2}}\Lambda(n)
e\left(\frac {ak^2n}{q}\right)
\end{eqnarray*}
Using (\ref{rep}), we get
\begin{eqnarray*}
\sum_{k<\sqrt{\frac{x}{u}}\atop (k,q)=1}\psi_u\left(\frac{x}{k^2},\chi\right)
&=& \frac {1}{\tau(\chi)}\sum_{a=1}^{q}\chi(a)T^\sharp\left(\frac{a}{q}\right)+
\frac {1}{\tau(\chi)}\sum_{a=1}^{q}\chi(a)T^\flat\left(\frac{a}{q}\right)\\
&=& \mathcal{O}\left(\left(q\sqrt q+\frac{x}{\sqrt q}+u^2\sqrt q\right)\log^2 x\right)\\
& &{}+\sum_{u<n\leqslant x}\sum_{u<m\leqslant x/n}\rho_u(m)\Lambda(n)\chi(mn)\\
&=& \mathcal{O}\left(\left(q\sqrt q+\frac{x}{\sqrt q}+u^2\sqrt q\right)\log^2 x\right)\\
& &{}+\sum_{u<m\leqslant x/u}\rho_u(m)\chi(m)\psi_u\left(\frac{x}{m},\chi\right).
\end{eqnarray*}
$\square$\\

The following lemma shows that the absolute value of the means of some suitable multiplicative functions does vary slowly.
\begin{lem}[\cite{GS1} corollary 3]\label{lem:absmean}
Let $f$ be a complex-valued multiplicative function with $|f(n)|\leqslant 1$. Then for $1\leqslant \omega\leqslant x/10$, we have $$\frac{1}{x}\left|\sum_{n\leqslant x}f(n)\right|-\frac{\omega}{x}\left|\sum_{n\leqslant x/\omega}f(n)\right|
\ll\left(\frac{\log(2\omega)}{\log x}\right)^{1-\frac{2}{\pi}}\log \left(\frac{\log x}{\log (2\omega)}\right)+\frac{\log\log x}{(\log x)^{2-\sqrt 3}},$$
The implied constant is absolute and computable.
\end{lem}

From the above lemma, we can deduce the following corollary.
\begin{cor}\label{cor:real}
Suppose $f$ is a real-valued multiplicative function with $|f(n)|\leqslant 1$ for all $n$. Then for $1\leqslant \omega\leqslant \sqrt x/2$, we have
$$\frac{1}{x}\sum_{n\leqslant x}f(n)-\frac{\omega}{x}\sum_{n\leqslant x/\omega}f(n)
\ll\left(\frac{\log(2\omega)}{\log x}\right)^{1-\frac{2}{\pi}}\log \left(\frac{\log x}{\log (2\omega)}\right)+\frac{\log\log x}{(\log x)^{2-\sqrt 3}}.$$
\end{cor}
Proof:
By Lemma~\ref{lem:absmean}, there exists an absolute constant $C_0>1$ such that
$$\left|\,\frac{1}{x}\Bigg|\sum_{n\leqslant x}f(n)\Bigg|-\frac{\omega}{x}\Bigg|\sum_{n\leqslant x/\omega}f(n)\Bigg|\,\right|
<C_0M(x,\omega),$$
where $$M(x,\omega)=\left(\frac{\log(2\omega)}{\log x}\right)^{1-\frac{2}{\pi}}\log \left(\frac{\log x}{\log (2\omega)}\right)+\frac{\log\log x}{(\log x)^{2-\sqrt 3}}.$$

If $\frac{1}{x}\left|\sum\limits_{n\leqslant x}f(n)\right|<2C_0M(x,\omega),$
then $\frac{\omega}{x}\left|\sum\limits_{n\leqslant x/\omega}f(n)\right|<3C_0M(x,\omega).$
Therefore
$$\left|\,\frac{1}{x}\sum_{n\leqslant x}f(n)-\frac{\omega}{x}\sum_{n\leqslant x/\omega}f(n)\,\right|
<5C_0M(x,\omega).$$
Otherwise, $$\frac{1}{x}\Bigg|\sum_{n\leqslant x}f(n)\Bigg|\geqslant 2C_0M(x,\omega).$$
Without loss of generality, we may assume
$$\frac{1}{x}\sum_{n\leqslant x}f(n)\geqslant 2C_0M(x,\omega).$$
Then $$\frac{\omega}{x}\sum_{n\leqslant x/\omega}f(n)>C_0M(x,\omega)\quad\text{or}\quad<-C_0M(x,\omega).$$
If $\frac{\omega}{x}\sum\limits_{n\leqslant x/\omega}f(n)>C_0M(x,\omega)$, we have
$$\left|\,\frac{1}{x}\sum_{n\leqslant x}f(n)-\frac{\omega}{x}\sum_{n\leqslant x/\omega}f(n)\,\right|=
\left|\,\frac{1}{x}\Bigg|\sum_{n\leqslant x}f(n)\Bigg|-\frac{\omega}{x}\Bigg|\sum_{n\leqslant x/\omega}f(n)\Bigg|\,\right|
<C_0M(x,\omega).$$
For the case $\frac{\omega}{x}\sum\limits_{n\leqslant x/\omega}f(n)<-C_0M(x,\omega)$, we can show that there exists a $x_0\in [\frac{x}{\omega},x]$ such that $\left|\sum\limits_{n\leqslant x_0}f(n)\right|\leqslant \frac{1}{2}$, since $\left(\sum\limits_{n\leqslant x}f(n)\right)\cdot \left(\sum\limits_{n\leqslant x/\omega}f(n)\right)<0$ and the real-valued function $f$ satisfies $|f(n)|\leqslant 1$ for all $n$.
Hence we get that
\begin{eqnarray}\label{ast}
  \left|\,\frac{1}{x}\bigg|\sum_{n\leqslant x}f(n)\bigg|-\frac{1}{x_0}\bigg|\sum_{n\leqslant x_0}f(n)\bigg|\,\right| &>& 2C_0M(x,\omega)-\frac{1}{2x_0}\\
    &>& 2C_0M(x,\omega)-\frac{\omega}{2x}.\nonumber
\end{eqnarray}
On the other hand, from Theorem~\ref{lem:absmean}, we deduce that
$$\left|\,\frac{1}{x}\bigg|\sum_{n\leqslant x}f(n)\bigg|-\frac{1}{x_0}\bigg|\sum_{n\leqslant x_0}f(n)\bigg|\,\right| < C_0M(x,x/x_0).$$
Since $x/x_0\leqslant \omega$, we have
\begin{equation}\label{aast}
\frac{\log x}{\log(2(x/x_0))}\geqslant\frac{\log x}{\log(2\omega)}.
\end{equation}
A basic observation is that the function $$\left(\frac{1}{t}\right)^{1-\frac{2}{\pi}}\log t$$ increases in $[1,\exp(\frac{1}{1-2/\pi})]$ and decreases in $[\exp(\frac{1}{1-2/\pi}),+\infty)$.
\begin{itemize}
  \item For $\frac{\log x}{\log(2\omega)}\geqslant\exp(\frac{1}{1-2/\pi})$, combining \ref{aast}, we have $M(x,x/x_0)\leqslant M(x,\omega)$. Hence
  $$\left|\,\frac{1}{x}\bigg|\sum_{n\leqslant x}f(n)\bigg|-\frac{1}{x_0}\bigg|\sum_{n\leqslant x_0}f(n)\bigg|\,\right| < C_0M(x,\omega).$$
  Furthermore, from $\frac{\log x}{\log(2\omega)}>\exp(\frac{1}{1-2/\pi})>15$, we have $2\omega<x^{1/15}$. It follows that for $x\geqslant 4$, we have
  $$\frac{\omega}{2x}<\frac{1}{4x^{14/15}}<\frac{\log\log x}{(\log x)^{2-\sqrt 3}}<C_0M(x,\omega),\qquad(\text{since}~C_0>1).$$
  Taking the inequality \ref{ast} into consideration, we get $$ \left|\,\frac{1}{x}\bigg|\sum_{n\leqslant x}f(n)\bigg|-\frac{1}{x_0}\bigg|\sum_{n\leqslant x_0}f(n)\bigg|\,\right| > C_0M(x,\omega),$$
  This leads to a contradiction!
  \item If $\frac{\log x}{\log (2\omega)}<\exp(\frac{1}{1-2/\pi})$, we get $2\leqslant \frac{\log x}{\log (2\omega)}<\exp(\frac{1}{1-2/\pi}),$ since $\omega\leqslant \sqrt x/2$.\\ Hence $\left(\frac{\log (2\omega)}{\log x}\right)^{1-\frac{2}{\pi}}\log\left(\frac{\log x}{\log(2\omega)}\right)\geqslant(\frac 12)^{1-\frac{2}{\pi}}\log 2$~.\\
  In such a case, the result is implied by the trivial estimate:
  \begin{equation*}
    \frac{1}{x}\sum_{n\leqslant x}f(n)-\frac{\omega}{x}\sum_{n\leqslant x/\omega}f(n) \ll\; 1\;\ll \;\left(\frac{\log (2\omega)}{\log x}\right)^{1-\frac{2}{\pi}}\log\left(\frac{\log x}
     {\log(2\omega)}\right).
  \end{equation*}
\end{itemize}
In conclusion, we always have
$$\frac{1}{x}\sum_{n\leqslant x}f(n)-\frac{\omega}{x}\sum_{n\leqslant x/\omega}f(n) \ll M(x,\omega).\qquad \square$$

\begin{lem}\label{LL}
Suppose $2\sqrt x<u^2<x$ and $\chi$ is a non-principle Dirichlet character modulo $q$, we have
\begin{eqnarray*}
  \sum_{u<m\leqslant x/u}\rho_u(m)\chi(m)\left(\frac{x}{m}-u\right) &=& \left(u\sum_{n\leqslant \frac{x}{u}}\lambda(n)\chi(n)\right)
  \left(L(1,\chi)\log {\frac{x}{eu^2}}+L'(1,\chi)\right) \\
   & & {}+\mathcal{O}\left(\left(u^2\log{\frac{x}{u^2}}\right)\sqrt q\log q+\frac{x}{u}\right)\\
   & & {}+\mathcal{O}\left(\epsilon(x,u)\cdot x\left(\log q+\left(\log {\frac{x}{u^2}}\right)^2\right)\right),
\end{eqnarray*}
where $$\epsilon(x,u)=\left(\left(\log{\frac{x}{u}}\right)^{\sqrt 3-2}+\left(\frac{\log {\frac{2x}{u^2}}}{\log{\frac{x}{u}}}\right)^{1-\frac{2}{\pi}}\right)\cdot \log\log x.$$
\end{lem}
Proof: By partial summation,
\begin{eqnarray*}
  \sum_{u<m\leqslant x/u}\rho_u(m)\chi(m)\left(\frac{x}{m}-u\right) &=& \sum_{u<t\leqslant \frac{x}{u}-1}\frac{x}{t(t+1)}
  \sum_{u<m\leqslant t}\rho_u(m)\chi(m) \\
   & & {}+\left(\frac{x}{[x/u]}-u\right)\sum_{u<m\leqslant \frac{x}{u}}\rho_u(m)\chi(m).
\end{eqnarray*}
From the property of $\lambda(n)$, we have
\begin{eqnarray*}
  \sum_{u<m\leqslant x/u}\rho_u(m)\chi(m) &=& \sum_{u<t^2\leqslant \frac{x}{u}\atop t>0}\chi(t^2)
-\sum_{u<m\leqslant x/u}\left(\sum_{d\mid m\atop d\leqslant u}\lambda(d)\right)\chi(m) \\
   &=&\mathcal{O}\left(\sqrt{\frac{x}{u}}\right)-\sum_{d\leqslant u}\lambda(d)\chi(d)\sum_{\frac{u}{d}<n\leqslant \frac{x}{u}}\chi(n)\\
   &=&\mathcal{O}\left(\sqrt{\frac{x}{u}}+u\sqrt q\log q\right),
\end{eqnarray*}
In the last step, we used the P\"{o}lya-Vinogradov estimate.\\
On the other hand, we have
\begin{eqnarray*}
  \sum_{u<m\leqslant t}\rho_u(m)\chi(m) &=& \sum_{u<m\leqslant t}\chi(m)\lambda(m)\sum_{d|m\atop d\leqslant \frac{m}{u}}\lambda(d)\\
    &=& \sum_{u<m\leqslant t}\rho_u(m)\chi(m)=\sum_{d\leqslant \frac{t}{u}}\chi(d)\sum_{u\leqslant n\leqslant \frac{t}{d}}\lambda(n)\chi(n).
\end{eqnarray*}
Combining the above identities, we have
\begin{eqnarray}\label{eq:sigma}
  {}\qquad\sum_{u<m\leqslant x/u}\rho_u(m)\chi(m)\left(\frac{x}{m}-u\right) &=& \sum_{u<t\leqslant \frac{x}{u}-1}\frac{x}{t(t+1)}
  \left(\sum_{d\leqslant \frac{t}{u}}\chi(d)\sum_{u\leqslant n\leqslant \frac{t}{d}}\lambda(n)\chi(n)\right)\\
  & &{}+\mathcal{O}\left(\frac{u^2}{x}\left(\sqrt{\frac{x}{u}}+u\sqrt q\log q\right)\right)\nonumber\\
   &:=& I_1-I_2+\mathcal{O}\left(\sqrt{\frac{x}{u}}+u\sqrt q\log q\right),\nonumber
\end{eqnarray}
where
$$I_1=\sum_{u<t\leqslant \frac{x}{u}-1}\frac{x}{t(t+1)}
  \left(\sum_{d\leqslant \frac{t}{u}}\chi(d)\sum_{n\leqslant \frac{t}{d}}\lambda(n)\chi(n)\right),$$

$$I_2=\sum_{u<t\leqslant \frac{x}{u}-1}\frac{x}{t(t+1)}
  \left(\sum_{d\leqslant \frac{t}{u}}\chi(d)\sum_{n<u}\lambda(n)\chi(n)\right).$$

\begin{eqnarray*}
    I_2 &=& x\sum_{n<u}\lambda(n)\chi(n)\left(\sum_{d\leqslant \frac{x/u-1}{u}}\chi(d)\sum_{du\leqslant t \leqslant \frac{x}{u}-1\atop u<t}\frac{1}{t(t+1)}\right) \\
      &=& x\sum_{n<u}\lambda(n)\chi(n)\left(\sum_{d\leqslant \frac{x/u-1}{u}}\chi(d)\left(\frac{1}{du}-\frac{1}{x/u}\right)+\mathcal{O}\left(\sum_{d\leqslant \frac{x/u-1}{u}}\frac{1}{(du)^2}\right)\right) \\
      &=& x\sum_{n<u}\lambda(n)\chi(n)\left(\frac{1}{u}\sum_{d\leqslant \frac{x/u-1}{u}}\frac{\chi(d)}{d}
          +\mathcal{O}\left(\frac{u\sqrt q\log q}{x}\right)+\mathcal{O}\left(\frac{1}{u^2}\right)\right)
\end{eqnarray*}
Noticing that
\begin{equation}\label{eq:appr}
\sum_{d\leqslant \frac{x/u-1}{u}}\frac{\chi(d)}{d}=L(1,\chi)+\mathcal{O}\left(\frac{u^2}{x}\sqrt q\log q\right),
\end{equation}
we have
\begin{equation}\label{eq:I2}
I_2=\frac{x}{u}L(1,\chi)\sum_{n<u}\lambda(n)\chi(n)+\mathcal{O}\left(u^2\sqrt q\log q+\frac{x}{u}\right).
\end{equation}
Applying Corollary~\ref{cor:real}, we have
\begin{equation}\label{eq:eqI2}
I_2=\left(u\sum_{n\leqslant \frac{x}{u}}\lambda(n)\chi(n)\right)\cdot L(1,\chi)+\mathcal{O}\left(u^2\sqrt q\log q+\frac{x}{u}+ \epsilon(x,u)x\log q\right).
\end{equation}
By Corollary ~\ref{cor:real}, for $du\leqslant t\leqslant \frac{x}{u}-1$, we have
\begin{equation}\label{eq:ide}
\frac{1}{t/d}\sum_{n\leqslant \frac{t}{d}}\lambda(n)\chi(n)=\frac{1}{x/u}\sum_{n\leqslant \frac{x}{u}}\lambda(n)\chi(n)
+\mathcal{O}\left(\epsilon(x,u)\right),
\end{equation}
where $$\epsilon(x,u)=\left(\left(\log{\frac{x}{u}}\right)^{\sqrt 3-2}+\left(\frac{\log {\frac{2x}{u^2}}}{\log{\frac{x}{u}}}\right)^{1-\frac{2}{\pi}}\right)\cdot \log\log x.$$

Due to the identity (\ref{eq:ide}), now we can handle the sum $I_1$.
\begin{eqnarray*}
  I_1 &=& \left(\frac{u}{x}\sum_{n\leqslant \frac{x}{u}}\lambda(n)\chi(n)\right)
  \left(\sum_{u<t\leqslant \frac{x}{u}-1}\frac{x}{t+1}\sum_{d\leqslant \frac{t}{u}}\frac{\chi(d)}{d}\right) \\
    & & {}+ \mathcal{O}\left(\epsilon(x,u)\left(\sum_{u<t\leqslant \frac{x}{u}-1}\frac{x}{t(t+1)}\sum_{d\leqslant \frac{t}{u}}\frac{t}{d}\right)\right)\\
    &=& \left(\frac{u}{x}\sum_{n\leqslant \frac{x}{u}}\lambda(n)\chi(n)\right)
  \left(\sum_{d\leqslant \frac{x/u-1}{u}}\frac{\chi(d)}{d}\sum_{du\leqslant t\leqslant \frac{x}{u}-1\atop u<t}\frac{x}{t+1}\right) \\
    & & {}+ \mathcal{O}\left(\epsilon(x,u)x\log^2\left(\frac{x}{u^2}\right)\right)\\
    &=& \left(\frac{u}{x}\sum_{n\leqslant \frac{x}{u}}\lambda(n)\chi(n)\right)
  \left(x\sum_{d\leqslant \frac{x/u-1}{u}}\frac{\chi(d)}{d}\left(\log\frac{x}{u^2}-\log d\right)+ \mathcal{O}\left(\sum_{d\leqslant \frac{x/u-1}{u}}\frac{x}{d^2u}\right)\right) \\
    & & {}+ \mathcal{O}\left(\epsilon(x,u)x\log^2\left(\frac{x}{u^2}\right)\right)\\
   &=& \left(u\sum_{n\leqslant \frac{x}{u}}\lambda(n)\chi(n)\right)
   \left(\log\frac{x}{u^2}\sum_{d\leqslant \frac{x/u-1}{u}}\frac{\chi(d)}{d}
   -\sum_{d\leqslant \frac{x/u-1}{u}}\frac{\chi(d)}{d}\log d\right) \\
    & & {}+\mathcal{O}\left(\frac{x}{u}+\epsilon(x,u)x\log^2\left(\frac{x}{u^2}\right)\right)\\
\end{eqnarray*}

Since
$$-\sum_{d\leqslant \frac{x/u-1}{u}}\frac{\chi(d)}{d}\log d
=L'(1,\chi)+\mathcal{O}\left(\frac{u^2\log\frac{x}{u^2}}{x}\sqrt q\log q\right),$$
combining (\ref{eq:appr}), we get
\begin{eqnarray}\label{eq:I1}
I_1 &=& \left(u\sum_{n\leqslant \frac{x}{u}}\lambda(n)\chi(n)\right)
   \left(\log\frac{x}{u^2}\cdot L(1,\chi)+L'(1,\chi)\right)\\
    & & {}+\mathcal{O}\left(\frac{x}{u}+\left(u^2\log\frac{x}{u^2}\right)
           \sqrt q\log q+\epsilon(x,u)x\log^2\left(\frac{x}{u^2}\right)\right).\nonumber
\end{eqnarray}

Finally, combing equations (\ref{eq:sigma}), (\ref{eq:I1}) and (\ref{eq:eqI2}), we get
\begin{eqnarray*}
\lefteqn{\sum_{u<m\leqslant x/u}\rho_u(m)\chi(m)\left(\frac{x}{m}-u\right)}\\
  &=& \left(u\sum_{n\leqslant \frac{x}{u}}\lambda(n)\chi(n)\right)
   \left(\log\frac{x}{eu^2}\cdot L(1,\chi)+L'(1,\chi)\right) \\
    & & {}+\mathcal{O}\left(\frac{x}{u}+\left(u^2\log\frac{x}{u^2}\right)
           \sqrt q\log q+\epsilon(x,u)x\left(\log q+\log^2\left(\frac{x}{u^2}\right)\right)\right).\quad\square
\end{eqnarray*}

\begin{lem}[see (22.109) in Chapter 22 of ~\cite{IK}]\label{lem:bound}
Suppose $\chi$ is a non-principle Dirichlet character modulo $q$ and $x\geqslant q$, we have
$$\sum_{n\leqslant x}\frac{\tau(n,\chi)}{n}=L(1,\chi)(\log x+\gamma)+L'(1,\chi)+\mathcal{O}\left(q^{1/4}x^{-1/2}\log x\right),$$
where $\tau(n,\chi)=\sum\limits_{d|n}\chi(d).$
\end{lem}

The following corollary shows that if $L(1,\chi)$ takes small value, then $\chi(p)$ takes negative value for most of small prime $p$.

\begin{cor}\label{cor:fapp}
Suppose $\chi$ is a quadratic Dirichlet character modulo $q$ and $x\geqslant q$. Then
$$\sum_{n\leqslant x}\Lambda(n)\chi(n)=-x+\mathcal{O}\left(\left(L(1,\chi)+q^{-1/4}\right)x\log^2 x+xe^{-c\sqrt{\log x}}+q\right),$$
where $c$ is some positive constant.
\end{cor}
Proof: By Lemma~\ref{lem:bound}, we have
\begin{equation}\label{eq:base}
\sum_{q<n\leqslant y}\frac{\tau(n,\chi)}{n}=L(1,\chi)\log {\frac{y}{q}}+\mathcal{O}\left(q^{-1/4}\log q\right).
\end{equation}
For a real character $\chi$, $\tau(n,\chi)\geqslant 0$. Then for $y>q$, we have
\begin{equation}\label{eq:bound}
\sum_{q<p\leqslant y}\frac{1+\chi(p)}{p}\log p\leqslant \log y\sum_{q<n\leqslant y}\frac{\tau(n,\chi)}{n}\ll \left(L(1,\chi)+q^{-1/4}\right)\log^2 y.
\end{equation}
Applying partial summation, we have
\begin{eqnarray}\label{eq:fapp}
  \sum_{q<p\leqslant x}(1+\chi(p))\log p &=& [x]\cdot\sum_{q<p\leqslant x}\frac{(1+\chi(p))\log p}{p}\\
    & & {}-\sum_{q<m\leqslant x-1}\sum_{q<p\leqslant m}\frac{(1+\chi(p))\log p}{p}\nonumber\\
    &\ll&\left(L(1,\chi)+q^{-1/4}\right)x\log^2 x.\nonumber
\end{eqnarray}
On the other hand,
\begin{eqnarray*}
  \sum_{n\leqslant x}\Lambda(n)\chi(n)&=& \sum_{n\leqslant x}\Lambda(n)(1+\chi(n))-\sum_{n\leqslant x}\Lambda(n)\\
    &=& \sum_{p\leqslant x}(1+\chi(p))\log p-x+\mathcal{O}\left(xe^{-c\sqrt{\log x}}\right)\\
    &=& -x+\sum_{q<p\leqslant x}(1+\chi(p))\log p+\mathcal{O}\left(q+xe^{-c\sqrt{\log x}}\right)
\end{eqnarray*}
Combining the estimate (\ref{eq:fapp}), we get the conclusion.$\qquad\square$

\begin{lem}[Proposition 4.5 in ~\cite{GS2}~ or Proposition 3 in ~\cite{GS1}~]\label{lem:multiplicative}
For any multiplicative function $f$ with $|f(p^k)|\leqslant 1$ for every prime power $p^k$, let
$$\Theta(f,x):=\prod_{p\leqslant x}\left(1+\frac{f(p)}{p}+\frac{f(p^2)}{p^2}+\cdots\right)\left(1-\frac{1}{p}\right),$$
and
$$s(f,x):=\sum_{p\leqslant x}\frac{|1-f(p)|}{p}.$$
For any $\varepsilon$ satisfying $1>\varepsilon\geqslant \frac{\log 2}{\log x}$, let $g$ be a completely multiplicative function satisfying the following condition:
\begin{equation*}
  g(p)=
  \begin{cases}
     1, & p\leqslant x^{\varepsilon}; \\
     f(p), & p> x^{\varepsilon}.
   \end{cases}
\end{equation*}
Then we have
$$\frac{1}{x}\sum_{n\leqslant x}f(n)=\Theta(f,x^\varepsilon)\frac{1}{x}\sum_{m\leqslant x}g(m)+\mathcal{O}(\varepsilon\exp(s(f,x))),$$
where the implied constant is absolute.
\end{lem}

Similarly as corollary ~\ref{cor:fapp}, we get the following result.
\begin{cor}\label{cor:fap}
Suppose $\chi$ is a quadratic Dirichlet character modulo $q$ and $x\geqslant q$. Then
$$\sum_{n\leqslant x}\chi(n)\lambda(n)=P(q)x+\mathcal{O}\left(\left(L(1,\chi)+q^{-1/4}\right)
x\log x+x\frac{\log^3 q}{\log x}\right),$$
where $$P(q)=\prod_{p\leqslant q}\left(1-\frac{1}{p}\right)\left(1+\frac{\chi(p)}{p}\right)^{-1}.$$
\end{cor}
Proof: Let the completely multiplicative function $f$ uniquely determined by the following condition:
\begin{equation*}
  f(p)=
  \begin{cases}
     \chi(p)\lambda(p)=-\chi(p), & p\leqslant q; \\
     1, & p> q.
   \end{cases}
\end{equation*}
Denote
$$E(x)=\{n\leqslant x\,|\,\text{ $n$ has at least a prime factor $p$ such that $p>q$ and $\lambda(p)\chi(p)\neq 1$~}\}.$$
Obviously, we have
\begin{eqnarray*}
  \left|\sum_{n\leqslant x}\chi(n)\lambda(n)-\sum_{n\leqslant x}f(n)\right| &\leqslant& 2|E(x)| \\
    &\ll& \sum_{q<p\leqslant x\atop \chi(p)\neq -1}\sum_{n\leqslant x\atop p\mid n}1 \\
    &\ll& x\sum_{q<p\leqslant x\atop \chi(p)\neq -1}\frac{1}{p}\\
    &\ll&  x\sum_{q<p\leqslant x}\frac{1+\chi(p)}{p}.
\end{eqnarray*}
Similarly as (\ref{eq:bound}), from (\ref{eq:base}) we get that
$$\sum_{q<p\leqslant x}\frac{1+\chi(p)}{p}\leqslant \sum_{q<n\leqslant x}\frac{\tau(n,\chi)}{n}\ll \left(L(1,\chi)+q^{-1/4}\right)
\log x.$$
Hence
\begin{equation}\label{eq:two}
\sum_{n\leqslant x}\chi(n)\lambda(n)=\sum_{n\leqslant x}f(n)+\mathcal{O}\left(\left(L(1,\chi)+q^{-1/4}\right)
x\log x\right).
\end{equation}
Now let the completely multiplicative function $g(n)\equiv 1$~.\\
Note that $g(p)=f(p)$ for $p>q$, by Lemma ~\ref{lem:multiplicative} we have
$$\frac{1}{x}\sum_{n\leqslant x}f(n)=\Theta\left(f,q\right)\frac{1}{x}\sum_{m\leqslant x}g(m)+\mathcal{O}\left(\frac{\log q}{\log x}\exp(s(f,x))\right),$$

where $$\Theta\left(f,q\right)=\prod_{p\leqslant q}\left(1-\frac{1}{p}\right)\left(1+\frac{\chi(p)}{p}\right)^{-1},$$
$$s(f,x)=\sum_{p\leqslant x}\frac{1-f(p)}{p}=\sum_{p\leqslant q}\frac{1+\chi(p)}{p}.$$

Therefore, $$\sum_{n\leqslant x}f(n)=\Theta\left(f,q\right)x+\mathcal{O}\left(x\frac{\log^3 q}{\log x}\right).$$

Combining (\ref{eq:two}), we have
$$\sum_{n\leqslant x}\chi(n)\lambda(n)=\Theta\left(f,q\right)x+\mathcal{O}\left(\left(L(1,\chi)+q^{-1/4}\right)
x\log x+x\frac{\log^3 q}{\log x}\right).\qquad \square$$

\textbf{The proof of Theorem ~\ref{main4}}:
Now we compute the sum
$$\sum_{u<m\leqslant x/u}\rho_u(m)\chi(m)\psi_u\left(\frac{x}{m},\chi\right)$$
 in two different ways under our assumption, which leads to the main result.\\
Let $T=\exp\left((\log q)^{8}\right)$~. We choose $x=T$ and $u=\sqrt{\frac{T}{q}}$~.
For sufficiently large $q$, we have $u\geqslant q$ and $u^2>2\sqrt x$.
From Lemma ~\ref{lem:psi} and Corollary ~\ref{cor:fapp}, it follows that
\begin{eqnarray*}
\lefteqn{\sum_{k<(qT)^{1/4}\atop (k,q)=1}\left(-\frac{T}{k^2}+\mathcal{O}\left(\sqrt{T/q}\right)\right)}\\
  & & {}+\mathcal{O}\left(\sum_{k<(qT)^{1/4}}\left(\left(L(1,\chi)+q^{-1/4}\right)\frac{T}{k^2}\log^2 T+\frac{T}{k^2}e^{-c\sqrt{\log (T/k^2)}}+q\right)\right)\\
 &=&-\sum_{\sqrt{\frac{T}{q}}<m\leqslant (qT)^{1/2}}\rho_u(m)\chi(m)\left(\frac{T}{m}-\sqrt{T/q}\right)+\mathcal{O}\left(\frac{T}{\sqrt q}\log^2 T\right) \\
    & & {}+\mathcal{O}\left(\sum_{\sqrt{\frac{T}{q}}<m\leqslant (qT)^{1/2}}d(m)\left(\left(L(1,\chi)+q^{-1/4}\right)\frac{T}{m}\log^2 T+\frac{T}{m}e^{-c\sqrt{\log (T/m)}}+q\right)\right).
\end{eqnarray*}
Since $\sum\limits_{m\leqslant y}d(m)\ll y\log y$, by partial summation we have
$$\sum_{\sqrt{\frac{T}{q}}<m\leqslant (qT)^{1/2}}\frac{d(m)}{m}\ll \log T\log q,$$
and
$$\sum_{\sqrt{\frac{T}{q}}<m\leqslant (qT)^{1/2}}\frac{d(m)}{m}e^{-c\sqrt{\log (T/m)}}\ll(\log T)^2e^{-c\sqrt{\frac{1}{2}\log (T/q)}}.$$
Hence
\begin{eqnarray}\label{identity}
T\sum_{k<(qT)^{1/4}\atop (k,q)=1}\frac{1}{k^2} &=& \sum_{\sqrt{\frac{T}{q}}<m\leqslant (qT)^{1/2}}\rho_u(m)\chi(m)\left(\frac{T}{m}-\sqrt{T/q}\right)\\
  & &{}+\mathcal{O}\left(\left(L(1,\chi)+q^{-1/4}\right)T(\log T)^3\log q\right).\nonumber
\end{eqnarray}
On the other hand, it follows from Lemma ~\ref{LL}~ and Corollary ~\ref{cor:fap}~that
\begin{eqnarray*}
  \lefteqn{\sum_{\sqrt{\frac{T}{q}}<m\leqslant (qT)^{1/2}}\rho_u(m)\chi(m)\left(\frac{T}{m}-\sqrt{T/q}\right)}\\
    &=& \left(P(q)T+\mathcal{O}\left(\left(L(1,\chi)+q^{-1/4}\right)T\log T+T\frac{\log^3 q}{\log T}\right)\right)\\
    & &{}\times\left(L'(1,\chi)+\mathcal{O}(L(1,\chi)\log q)\right) \\
    & &{}+\mathcal{O}\left(\frac{T}{\sqrt q}\log^2 q+T\log^2 q\left((\log T)^{\sqrt 3-2}+\left(\frac{\log q} {\log T}\right)^{1-\frac{2}{\pi}}\right)\cdot \log\log T\right).
\end{eqnarray*}
Note that $\log T=(\log q)^8$ and $P(q)\gg (\log q)^{-2}$, under the assumption of $L(1,\chi)\ll (\log q)^{-26}$, we have
\begin{equation}\label{indrect}
\sum_{\sqrt{\frac{T}{q}}<m\leqslant (qT)^{1/2}}\rho_u(m)\chi(m)\left(\frac{T}{m}-\sqrt{T/q}\right)= P(q)L'(1,\chi)T+\mathcal{O}\left(T(\log q)^{-1/10}\right).
\end{equation}
Besides,
\begin{equation}\label{coefficient}
  \sum_{k<(qT)^{1/4}\atop (k,q)=1}\frac{1}{k^2}=\frac{\pi^2}{6}\prod_{p|q}\left(1-\frac{1}{p^2}\right)+\mathcal{O}\left((qT)^{-1/4}\right).
\end{equation}
Inserting (\ref{indrect}) and (\ref{coefficient}) into (\ref{identity}) gives
$$P(q)L'(1,\chi)=\frac{\pi^2}{6}\prod_{p|q}\left(1-\frac{1}{p^2}\right)+\mathcal{O}\left((\log q)^{-1/10}\right).$$
It follows that
$$L'(1,\chi)=\frac{\pi^2}{6}\prod_{p|q}\left(1+\frac{1}{p}\right)\prod_{p\leqslant q\atop \chi(p)=1}\left(1+\frac{1}{p}\right)\left(1-\frac{1}{p}\right)^{-1}\Big(1+\mathcal{O}\left((\log q)^{-1/10}\right)\Big).$$
$\square$

E-mail address: lhyan@amss.ac.cn\par
Academy of Mathematics and System Science, Chinese Academy of Sciences, Beijing,\par
100190, P.R.China

\end{document}